\documentclass[11pt,a4paper]{amsart}
\usepackage{amsmath,amssymb,amsthm}
\usepackage{geometry}
\newtheorem{theorem}{Theorem}[section]
\newtheorem{proposition}[theorem]{Proposition}
\newtheorem{lemma}[theorem]{Lemma}
\theoremstyle{definition}
\newtheorem{definition}[theorem]{Definition}
\theoremstyle{remark}
\newtheorem{remark}[theorem]{Remark}

\title[Global Euler Cocycle in the BSS Complex of $SL(4,\mathbb{R})$]{Global Euler Cocycle in the Bott--Shulman--Stasheff Complex of $SL(4,\mathbb{R})$}
\author{Naoya Suzuki}
\date{November 29, 2026}
\subjclass[2020]{57R20, 55R40, 22E41, 58A10}
\keywords{Bott--Shulman--Stasheff complex, Euler class, simplicial de Rham cohomology, $SL(4,\mathbb{R})$ polar decomposition, Wigner rotation}

\begin{document}

\begin{abstract}
In a previous work \cite{Suzuki2016}, the author exhibited an explicit cocycle in the Bott--Shulman--Stasheff (BSS) simplicial de Rham complex $\Omega^*(NSO(4))$ representing the universal Euler class of $ESO(4)\to BSO(4)$. Due to the simplicial affine linearity of the canonical universal connection on compact orthogonal groups, that cocycle terminates at simplicial degree $2$, consisting solely of $(1,3)$- and $(2,2)$-forms. In this paper, we construct a smooth, globally defined $4$-cocycle
\[
E = E_{1,3} + E_{2,2} + E_{3,1} + E_{4,0} \in \Omega_{\mathrm{tot}}^4(NSL(4,\mathbb{R}))
\]
extending the $SO(4)$ Euler cocycle to the non-compact Lie group $G = SL(4,\mathbb{R})$. Because the polar orthogonal projection $u \colon SL(4,\mathbb{R}) \to SO(4)$ fails to be a group homomorphism, higher simplicial components $E_{3,1} \in \Omega^1(G^3)$ and a global smooth group $4$-cochain $E_{4,0} \in C^\infty(G^4)$ necessarily appear. By exploiting the global diffeomorphism of the matrix logarithm on positive-definite symmetric matrices and the right $SO(4)$-equivariance of polar decomposition, we define a family of nested hypercube associators $W_p \colon [0,1]^{p-1}\times G^p \to SO(4)$ whose boundary faces precisely govern both simplicial face contractions and multiplicative splittings. Pulling back the $SO(4)$ Euler cocycle along $W_p$ and integrating over the hypercubes yields a closed total $4$-cocycle whose restriction to $NSO(4)$ recovers the standard Euler cocycle and whose cohomology class represents the universal real Euler class in $H^4(BSL(4,\mathbb{R});\mathbb{R})$. Furthermore, we demonstrate by an explicit rational pure-boost $4$-tuple in $(SL(4,\mathbb{R})\setminus SO(4))^4$ that $E_{4,0}$ is non-zero as a smooth group $4$-cocycle.
\end{abstract}

\maketitle

\section{Introduction}

For any Lie group $G$, the nerve $NG = \{NG(p)\}_{p\ge 0}$ is a simplicial manifold with $NG(p) = G^p$, and its fat realization $\|NG\|$ models the classifying space $BG$ \cite{Segal1968}. Bott, Shulman, and Stasheff \cite{BSS1976} showed that the cohomology of the total complex $\Omega_{\mathrm{tot}}^*(NG)$ of the simplicial de Rham double complex $\Omega^{p,q}(NG) = \Omega^q(G^p)$ is canonically isomorphic to the real cohomology $H^*(BG;\mathbb{R})$ \cite{Bott1973,BSS1976,Dupont1978,MostowPerchick1973}. Using Dupont's simplicial differential forms $A^{*,*}(NG)$ \cite{Dupont1976}, the author previously constructed explicit representatives in $\Omega_{\mathrm{tot}}^*(NG)$ for the Chern characters \cite{Suzuki2015} and for the Euler classes of $ESO(2m)\to BSO(2m)$ \cite{Suzuki2016}.

For $K = SO(4)$, the universal Euler cocycle obtained in \cite[Theorem 3.2]{Suzuki2016} takes a remarkably compact form: it is supported entirely in bidegrees $(1,3)$ and $(2,2)$,
\[
\omega^{SO(4)} = \omega_{1,3} + \omega_{2,2} \in \Omega^3(SO(4)) \oplus \Omega^2(SO(4)^2),
\]
satisfying $d\omega_{1,3} = 0$, $d'\omega_{1,3} + d\omega_{2,2} = 0$, and $d'\omega_{2,2} = 0$. This truncation occurs because the canonical universal connection $\theta = \sum_{i=0}^p t_i \theta_i$ on $ESO(4)\to BSO(4)$ produces curvature components whose Pfaffian integrals vanish beyond $p=2$.

A natural question is how to extend this explicit Euler cocycle from $SO(4)$ to the ambient non-compact Lie group $G = SL(4,\mathbb{R})$. Since the inclusion $\iota \colon SO(4) \hookrightarrow SL(4,\mathbb{R})$ is a homotopy equivalence, the induced map $H_{\mathrm{tot}}^4(NSL(4,\mathbb{R})) \to H_{\mathrm{tot}}^4(NSO(4))$ is an isomorphism. However, constructing an explicit global representative on $NSL(4,\mathbb{R})$ encounters two fundamental obstructions:
\begin{enumerate}
\item \emph{Failure of multiplicativity (Wigner rotation defect)}: The polar decomposition $h = u(h)p(h)$ provides a global smooth retraction $u \colon SL(4,\mathbb{R}) \to SO(4)$, where $u(h) = h(h^T h)^{-1/2}$. Simply pulling back $\omega_{1,3}$ and $\omega_{2,2}$ via $u$ fails to define a cocycle because $u(h_1 h_2) \neq u(h_1)u(h_2)$ in general. Restoring the $(2,3)$-cocycle equation requires adding a homotopy transgression term to $E_{2,2}$, which in turn creates a non-zero horizontal coboundary $d'E_{2,2} \neq 0$ in bidegree $(3,2)$. Canceling this cascade inevitably forces the appearance of a $1$-form $E_{3,1} \in \Omega^1(G^3)$ and a real-valued smooth function $E_{4,0} \in C^\infty(G^4)$.
\item \emph{Global smoothness vs.\ local logarithm}: Naive attempts to construct homotopy primitives using the Lie algebra logarithm $\log h$ on $SL(4,\mathbb{R})$ are obstructed by the branch cut along matrices with negative real eigenvalues, restricting such formulas to a small neighborhood of the identity.
\end{enumerate}

In this paper, we overcome both obstructions globally on $G = SL(4,\mathbb{R})$. The key geometric observation is twofold: first, while $\log h$ is only locally defined on $SL(4,\mathbb{R})$, the symmetric factor $p(h) = (h^T h)^{1/2}$ lies in the space $\mathcal{P}_1(4)$ of positive-definite symmetric matrices of determinant $1$, on which the principal matrix logarithm $\log \colon \mathcal{P}_1(4) \to \mathfrak{p}$ is a global $C^\infty$ diffeomorphism. Second, although $u \colon G \to SO(4)$ is not a homomorphism, it satisfies strict right $SO(4)$-equivariance, $u(gk) = u(g)k$ for all $g \in G$ and $k \in SO(4)$.

Using the global geodesic interpolation $\eta(\tau; h) = u(h)\exp(\tau \log p(h))$ ($\tau \in [0,1]$), we define a family of nested hypercube associators
\[
W_p \colon [0,1]^{p-1} \times G^p \longrightarrow SO(4) \qquad (p \ge 1).
\]
On the boundary facets of the $(p-1)$-cube $[0,1]^{p-1}$, the maps $W_p$ satisfy exact algebraic identities: setting $\tau_k = 1$ contracts the $k$-th adjacent product via the simplicial face operator $\epsilon_k$, whereas setting $\tau_k = 0$ splits $W_p$ into a product $W_k \cdot W_{p-k}$ in $SO(4)$ (Proposition~\ref{prop:boundary}). Pulling back $\omega_{1,3}$ and $\omega_{2,2}$ along $W_p$ and integrating over the fiber hypercubes $[0,1]^{p-1}$ yields explicit closed formulas for $E_{1,3}, E_{2,2}, E_{3,1}, E_{4,0}$ (Section~\ref{sec:explicit}), which satisfy $(d'+d'')E = 0$ globally on $NSL(4,\mathbb{R})$ (Theorem~\ref{thm:main}).

The remainder of the paper is organized as follows. Section~\ref{sec:prelim} establishes the conventions for the BSS double complex, fiber integration, and the baseline $SO(4)$ Euler cocycle. Section~\ref{sec:associators} introduces the nested hypercube associators $W_p$ and proves their boundary factorization formulas. Section~\ref{sec:explicit} gives the explicit integral formulas for $E_{1,3}$, $E_{2,2}$, $E_{3,1}$, $E_{4,0}$, discusses vanishing on central elements, and presents an explicit rational pure-boost $4$-tuple in $(SL(4,\mathbb{R})\setminus SO(4))^4$ on which $E_{4,0}$ evaluates to $\frac{1}{8\pi^2}(\arctan\frac{9}{25})^2 \neq 0$ (Proposition~\ref{prop:nonvanishing}). Section~\ref{sec:proof} proves Theorem~\ref{thm:main} across all bidegrees $(1,4)$ through $(5,0)$. Section~\ref{sec:remarks} remarks on Brylinski's truncated local complex \cite{Brylinski2000} and the loop Lie algebra $L\mathfrak{sl}_4(\mathbb{R})$.

\section{Preliminaries: The BSS Double Complex and the $SO(4)$ Euler Cocycle}\label{sec:prelim}

\subsection{Conventions for the Simplicial de Rham Complex}
Let $G = SL(4,\mathbb{R})$. For each $p \ge 0$, let $NG(p) = G^p$ (with $NG(0) = \{*\}$) be the $p$-th level of the nerve of $G$, equipped with the face operators $\epsilon_i \colon G^p \to G^{p-1}$ ($0 \le i \le p$) defined by \cite{BSS1976,Dupont1978,Suzuki2016}
\[
\epsilon_i(h_1, \dots, h_p) =
\begin{cases}
(h_2, \dots, h_p) & (i = 0), \\
(h_1, \dots, h_i h_{i+1}, \dots, h_p) & (1 \le i \le p-1), \\
(h_1, \dots, h_{p-1}) & (i = p).
\end{cases}
\]
These satisfy the standard simplicial identities $\epsilon_j \epsilon_i = \epsilon_{i-1} \epsilon_j$ for $0 \le j < i \le p$.

\begin{definition}\label{def:bss}
The Bott--Shulman--Stasheff (BSS) double complex $(\Omega^{*,*}(NG), d', d'')$ is given by $\Omega^{p,q}(NG) := \Omega^q(G^p)$ ($p \ge 1, q \ge 0$; $\Omega^{0,0}(NG) = \mathbb{R}$, $\Omega^{0,q}(NG) = 0$ for $q > 0$) with horizontal and vertical differentials
\begin{align*}
d' &:= \sum_{i=0}^p (-1)^i \epsilon_i^* \colon \Omega^q(G^{p-1}) \longrightarrow \Omega^q(G^p), \\
d'' &:= (-1)^p d_{G^p} \colon \Omega^q(G^p) \longrightarrow \Omega^{q+1}(G^p),
\end{align*}
where $d_{G^p}$ denotes the ordinary exterior derivative on $G^p$.
\end{definition}

One checks immediately that $(d')^2 = 0$, $(d'')^2 = 0$, and $d'd'' + d''d' = 0$. By the fundamental theorem of Bott--Shulman--Stasheff \cite{BSS1976,Dupont1978,MostowPerchick1973}, the cohomology of the total complex $(\Omega_{\mathrm{tot}}^*(NG), D = d'+d'')$ is canonically isomorphic to $H^*(BG;\mathbb{R})$.

For a total $4$-cochain $E = E_{1,3} + E_{2,2} + E_{3,1} + E_{4,0}$ with $E_{p,4-p} \in \Omega^{4-p}(G^p)$, the total cocycle condition $(d'+d'')E = 0$ decomposes by bidegree into the system of five equations:
\begin{align*}
(1)&\ \text{Bidegree } (1,4): & d_G E_{1,3} &= 0, \\
(2)&\ \text{Bidegree } (2,3): & d'E_{1,3} + d_{G^2} E_{2,2} &= 0, \\
(3)&\ \text{Bidegree } (3,2): & d'E_{2,2} - d_{G^3} E_{3,1} &= 0, \\
(4)&\ \text{Bidegree } (4,1): & d'E_{3,1} + d_{G^4} E_{4,0} &= 0, \\
(5)&\ \text{Bidegree } (5,0): & d'E_{4,0} &= 0.
\end{align*}

\subsection{Fiber Integration Convention and Stokes' Formula}
Let $I^r = [0,1]^r$ with coordinates $(\tau_1, \dots, \tau_r)$, and let $X$ be a smooth manifold. For a differential form $\alpha \in \Omega^{r+q}(I^r \times X)$, we define the fiber integral $\int_{I^r} \alpha \in \Omega^q(X)$ by contracting $d\tau_1 \wedge \cdots \wedge d\tau_r$ from the left: if
\[
\alpha = d\tau_1 \wedge \cdots \wedge d\tau_r \wedge \beta(\tau, x) + (\text{terms of degree } < r \text{ in } d\tau_1, \dots, d\tau_r),
\]
then
\[
\int_{I^r} \alpha := \int_{I^r} \beta(\tau, x)\, d\tau_1 \cdots d\tau_r \in \Omega^q(X).
\]

\begin{lemma}[Fiberwise Stokes' Formula]\label{lem:stokes}
Under the above convention, for any $\omega \in \Omega^{r+q-1}(I^r \times X)$,
\[
\int_{I^r} d_{I^r \times X}\omega = \sum_{k=1}^r (-1)^{k-1}\left(\int_{I^{r-1}} \omega|_{\tau_k=1} - \int_{I^{r-1}} \omega|_{\tau_k=0}\right) + (-1)^r d_X \int_{I^r} \omega,
\]
where $\int_{I^{r-1}}$ denotes integration over the remaining coordinates $(\tau_1, \dots, \hat{\tau}_k, \dots, \tau_r)$ (and evaluates to the restriction itself when $r = 1$).
\end{lemma}

\begin{proof}
Write
\[
\omega = d\tau_1 \wedge \cdots \wedge d\tau_r \wedge \gamma + \sum_{k=1}^r (-1)^{k-1} d\tau_1 \wedge \cdots \wedge \widehat{d\tau_k} \wedge \cdots \wedge d\tau_r \wedge \sigma_k + \rho,
\]
where $\gamma, \sigma_k$ contain no $d\tau_j$ and $\rho$ has degree at most $r-2$ in $d\tau_j$. Taking $d_{I^r \times X}\omega = (d_{I^r} + d_X)\omega$, the component proportional to $d\tau_1 \wedge \cdots \wedge d\tau_r$ is
\[
d\tau_1 \wedge \cdots \wedge d\tau_r \wedge \left((-1)^r d_X \gamma + \sum_{k=1}^r (-1)^{k-1} \frac{\partial \sigma_k}{\partial \tau_k}\right).
\]
Integrating over $I^r$ and applying the fundamental theorem of calculus in each $\tau_k$ yields the stated identity.
\end{proof}

\subsection{The Baseline Euler Cocycle on $SO(4)$}
Let $K = SO(4)$. Following \cite[Theorem 3.2]{Suzuki2016}, the universal Euler class of $ESO(4)\to BSO(4)$ in $\Omega_{\mathrm{tot}}^4(NSO(4))$ is represented by $\omega^{SO(4)} = \omega_{1,3} + \omega_{2,2}$, where
\begin{align*}
\omega_{1,3}(u) &:= \frac{1}{96\pi^2} \sum_{\sigma \in \mathfrak{S}_4} \operatorname{sgn}(\sigma) (u^{-1}du)_{\sigma(1)\sigma(2)} \wedge (u^{-1}du \wedge u^{-1}du)_{\sigma(3)\sigma(4)} \in \Omega^3(SO(4)), \\
\omega_{2,2}(u_1, u_2) &:= -\frac{1}{32\pi^2} \sum_{\sigma \in \mathfrak{S}_4} \operatorname{sgn}(\sigma) (u_1^{-1}du_1)_{\sigma(1)\sigma(2)} \wedge (du_2 u_2^{-1})_{\sigma(3)\sigma(4)} \in \Omega^2(SO(4)^2).
\end{align*}

\begin{lemma}[{\cite[Theorem 3.2]{Suzuki2016}}]\label{lem:so4cocycle}
On $K = SO(4)$, the forms $\omega_{1,3}$ and $\omega_{2,2}$ satisfy:
\begin{enumerate}
\item $d_K \omega_{1,3} = 0$.
\item $d'_K \omega_{1,3} + d_K \omega_{2,2} = 0$, i.e., for any smooth maps $k_1, k_2 \colon M \to SO(4)$ on a manifold $M$,
\[
(k_1 k_2)^* \omega_{1,3} = k_1^* \omega_{1,3} + k_2^* \omega_{1,3} + d_M\big((k_1, k_2)^* \omega_{2,2}\big).
\]
\item $d'_K \omega_{2,2} = 0$, i.e., for any three smooth maps $k_1, k_2, k_3 \colon M \to SO(4)$,
\[
(k_2, k_3)^* \omega_{2,2} - (k_1 k_2, k_3)^* \omega_{2,2} + (k_1, k_2 k_3)^* \omega_{2,2} - (k_1, k_2)^* \omega_{2,2} = 0.
\]
\end{enumerate}
\end{lemma}

\section{Nested Hypercube Associators and Their Boundary Formulas}\label{sec:associators}

\subsection{Polar Decomposition and Global Geodesic Interpolation}
Let $\mathcal{P}_1(4) := \{P \in \mathrm{Sym}_4^+(\mathbb{R}) \mid \det P = 1\}$ denote the symmetric space of $4\times 4$ real positive-definite symmetric matrices of determinant $1$, and let $\mathfrak{p} := \{Y \in \mathrm{Sym}_4(\mathbb{R}) \mid \operatorname{tr}(Y) = 0\}$. Since every $P \in \mathcal{P}_1(4)$ has strictly positive eigenvalues, the matrix square root $P \mapsto P^{1/2}$ and the principal matrix logarithm $\log \colon \mathcal{P}_1(4) \to \mathfrak{p}$ are global $C^\infty$ diffeomorphisms.

We define the global smooth maps $p \colon G \to \mathcal{P}_1(4)$ and $u \colon G \to SO(4)$ on $G = SL(4,\mathbb{R})$ by
\[
p(h) := (h^T h)^{1/2} \in \mathcal{P}_1(4), \qquad u(h) := h\, p(h)^{-1} \in SO(4) \qquad (h \in G).
\]
Then $h = u(h)p(h)$ is the global smooth polar decomposition of $G$.

\begin{lemma}[Bilateral $SO(4)$-Equivariance]\label{lem:equiv}
For all $g \in SL(4,\mathbb{R})$ and $k \in SO(4)$,
\[
u(gk) = u(g)k \qquad \text{and} \qquad u(kg) = k\, u(g).
\]
\end{lemma}

\begin{proof}
For right multiplication, $(gk)^T(gk) = k^T g^T g k$, so uniqueness of the positive-definite square root gives $p(gk) = k^T p(g) k$. Hence $u(gk) = gk(k^T p(g)^{-1} k) = g p(g)^{-1} k = u(g)k$. For left multiplication, $(kg)^T(kg) = g^T k^T k g = g^T g$, so $p(kg) = p(g)$ and $u(kg) = k g p(g)^{-1} = k u(g)$.
\end{proof}

\begin{definition}\label{def:eta}
Define the global smooth interpolation map $\eta \colon [0,1] \times G \to G$ by
\[
\eta(\tau; h) := u(h)\exp(\tau \log p(h)).
\]
Clearly $\eta$ is $C^\infty$ on $[0,1] \times G$ and satisfies:
\begin{itemize}
\item $\eta(0; h) = u(h) \in SO(4)$,
\item $\eta(1; h) = u(h)p(h) = h$,
\item $\eta(\tau; k) = k$ for all $\tau \in [0,1]$ and $k \in SO(4)$ (since $p(k) = I_4$ and $\log I_4 = 0$).
\end{itemize}
\end{definition}

\subsection{Nested Hypercube Associators $W_p$}
To interpolate the failure of multiplicativity across $p$-tuples $(h_1, \dots, h_p) \in G^p$, we assign one interpolation parameter $\tau_k \in [0,1]$ to each of the $p-1$ multiplicative gaps between $h_k$ and $h_{k+1}$.

\begin{definition}\label{def:Wp}
For each integer $p \ge 1$, define $W_p \colon [0,1]^{p-1} \times G^p \to SO(4)$ recursively by
\begin{align*}
W_1(h_1) &:= u(h_1), \\
W_2(\tau_1; h_1, h_2) &:= u\big(h_1 \cdot \eta(\tau_1; h_2)\big), \\
W_p(\tau_1, \dots, \tau_{p-1}; h_1, \dots, h_p) &:= u\big(h_1 \cdot \eta(\tau_1; h_2 \cdot \eta(\tau_2; \cdots h_{p-1} \cdot \eta(\tau_{p-1}; h_p) \cdots ))\big).
\end{align*}
More formally, setting $X_1(h_p) := h_p$ and $X_m(\tau_{p-m+1}, \dots, \tau_{p-1}; h_{p-m+1}, \dots, h_p) := h_{p-m+1} \cdot \eta(\tau_{p-m+1}; X_{m-1})$ for $2 \le m \le p$, we have $W_p = u \circ X_p$.
\end{definition}

Since $u$ and $\eta$ are globally $C^\infty$, each $W_p \colon [0,1]^{p-1} \times G^p \to SO(4)$ is a globally smooth map. Moreover, if $(h_1, \dots, h_p) \in SO(4)^p$, then $\eta(\tau_k; \cdot)$ acts as the identity and $W_p(\tau_1, \dots, \tau_{p-1}; h_1, \dots, h_p) \equiv h_1 h_2 \cdots h_p$ independently of $(\tau_1, \dots, \tau_{p-1})$.

\begin{proposition}[Boundary Factorization of $W_p$]\label{prop:boundary}
For any $p \ge 2$ and any $1 \le k \le p-1$, the restriction of $W_p$ to the boundary facets $\tau_k = 1$ and $\tau_k = 0$ of $[0,1]^{p-1}$ satisfies:
\begin{enumerate}
\item \emph{Contraction at $\tau_k = 1$}:
\[
W_p(\tau_1, \dots, 1_k, \dots, \tau_{p-1}; h_1, \dots, h_p) = W_{p-1}(\tau_1, \dots, \hat{\tau}_k, \dots, \tau_{p-1}; \epsilon_k(h_1, \dots, h_p)).
\]
\item \emph{Splitting at $\tau_k = 0$}:
\[
W_p(\tau_1, \dots, 0_k, \dots, \tau_{p-1}; h_1, \dots, h_p) = W_k(\tau_1, \dots, \tau_{k-1}; h_1, \dots, h_k) \cdot W_{p-k}(\tau_{k+1}, \dots, \tau_{p-1}; h_{k+1}, \dots, h_p).
\]
\end{enumerate}
\end{proposition}

\begin{proof}
Set $Y := X_{p-k}(\tau_{k+1}, \dots, \tau_{p-1}; h_{k+1}, \dots, h_p)$, so that $u(Y) = W_{p-k}(\tau_{k+1}, \dots, \tau_{p-1}; h_{k+1}, \dots, h_p)$.

For (1), at $\tau_k = 1$ we have $h_k \cdot \eta(1; Y) = h_k Y = (h_k h_{k+1}) \cdot \eta(\tau_{k+1}; X_{p-k-1})$ (or simply $h_{p-1} h_p$ when $k = p-1$). Thus the nested expression collapses the adjacent pair $(h_k, h_{k+1})$ into their product $h_k h_{k+1}$, yielding $W_{p-1}$ evaluated on $\epsilon_k(h_1, \dots, h_p)$.

For (2), at $\tau_k = 0$ we have $\eta(0; Y) = u(Y) = W_{p-k} \in SO(4)$. We claim by downward induction on $j = k, k-1, \dots, 1$ that
\[
X_{p-j+1}(\tau_j, \dots, 0_k, \dots, \tau_{p-1}; h_j, \dots, h_p) = \widetilde{X}_{k-j+1}(\tau_j, \dots, \tau_{k-1}; h_j, \dots, h_k) \cdot W_{p-k},
\]
where $\widetilde{X}$ denotes the nested expression for the $k$-tuple $(h_1, \dots, h_k)$. For $j = k$ the left-hand side is $h_k \cdot u(Y) = \widetilde{X}_1(h_k) \cdot W_{p-k}$. Suppose the claim holds for $j+1$. By Lemma~\ref{lem:equiv}, right multiplication by $K_0 := W_{p-k} \in SO(4)$ satisfies $u(Z K_0) = u(Z) K_0$ and $p(Z K_0) = K_0^T p(Z) K_0$, whence
\[
\eta(\tau_j; Z K_0) = u(Z) K_0 \exp\big(\tau_j K_0^T (\log p(Z)) K_0\big) = u(Z) \exp(\tau_j \log p(Z)) K_0 = \eta(\tau_j; Z) K_0.
\]
Multiplying on the left by $h_j$ advances the induction to $j$. At $j = 1$, applying $u$ and Lemma~\ref{lem:equiv} gives
\[
W_p|_{\tau_k = 0} = u(\widetilde{X}_k \cdot W_{p-k}) = u(\widetilde{X}_k) \cdot W_{p-k} = W_k \cdot W_{p-k}. \qedhere
\]
\end{proof}

\section{Explicit Formulas for the Global Euler Cocycle on $NSL(4,\mathbb{R})$}\label{sec:explicit}

Using the baseline forms $\omega_{1,3} \in \Omega^3(SO(4))$ and $\omega_{2,2} \in \Omega^2(SO(4)^2)$ from Section~\ref{sec:prelim} and the nested hypercube associators $W_p$ from Definition~\ref{def:Wp}, we define the components $E_{p,4-p} \in \Omega^{4-p}(G^p)$ ($p = 1,2,3,4$) of our total $4$-cochain $E$ on $NSL(4,\mathbb{R})$.

\begin{definition}\label{def:E}
Define $E = E_{1,3} + E_{2,2} + E_{3,1} + E_{4,0} \in \Omega_{\mathrm{tot}}^4(NSL(4,\mathbb{R}))$ by:
\begin{enumerate}
\item \emph{Bidegree $(1,3)$}: $E_{1,3} \in \Omega^3(G)$,
\[
E_{1,3}(h_1) := W_1(h_1)^* \omega_{1,3} = u(h_1)^* \omega_{1,3}.
\]
\item \emph{Bidegree $(2,2)$}: $E_{2,2} \in \Omega^2(G^2)$,
\[
E_{2,2}(h_1, h_2) := (W_1(h_1), W_1(h_2))^* \omega_{2,2} + \int_0^1 W_2(\tau_1; h_1, h_2)^* \omega_{1,3}.
\]
\item \emph{Bidegree $(3,1)$}: $E_{3,1} \in \Omega^1(G^3)$,
\begin{align*}
E_{3,1}(h_1, h_2, h_3) &:= \int_{[0,1]^2} W_3(\tau_1, \tau_2; h_1, h_2, h_3)^* \omega_{1,3} \\
&\quad - \int_0^1 (W_2(\tau_1; h_1, h_2), W_1(h_3))^* \omega_{2,2} + \int_0^1 (W_1(h_1), W_2(\tau_2; h_2, h_3))^* \omega_{2,2}.
\end{align*}
\item \emph{Bidegree $(4,0)$}: $E_{4,0} \in C^\infty(G^4)$ (a globally smooth real-valued function),
\begin{align*}
E_{4,0}(h_1, h_2, h_3, h_4) &:= \int_{[0,1]^3} W_4(\tau_1, \tau_2, \tau_3; h_1, h_2, h_3, h_4)^* \omega_{1,3} \\
&\quad + \int_{[0,1]^2} (W_3(\tau_1, \tau_2; h_1, h_2, h_3), W_1(h_4))^* \omega_{2,2} \\
&\quad - \int_{[0,1]^2} (W_2(\tau_1; h_1, h_2), W_2(\tau_3; h_3, h_4))^* \omega_{2,2} \\
&\quad + \int_{[0,1]^2} (W_1(h_1), W_3(\tau_2, \tau_3; h_2, h_3, h_4))^* \omega_{2,2}.
\end{align*}
\end{enumerate}
\end{definition}

\begin{remark}[Vanishing on $SO(4)^4$ and Central Elements]\label{rem:vanishing}
Because $\det(-I_4) = (-1)^4 = 1$ and $(-I_4)^T(-I_4) = I_4$, the central involution $-I_4$ lies in the maximal compact subgroup $SO(4) \subset SL(4,\mathbb{R})$ (with $p(-I_4) = I_4$ and $u(-I_4) = -I_4$). More generally, whenever $(h_1, h_2, h_3, h_4) \in SO(4)^4$---and in particular for $(h_1, h_2, h_3, h_4) = (I_4, I_4, I_4, I_4)$ or $(-I_4, -I_4, -I_4, -I_4)$---each map $W_p(\tau_1, \dots, \tau_{p-1}; h_1, \dots, h_p)$ is identically equal to $h_1 \cdots h_p$ independent of $(\tau_1, \dots, \tau_{p-1})$. Consequently, all fiber integrals of positive-degree forms vanish identically, giving
\[
E_{4,0}(I_4, I_4, I_4, I_4) = E_{4,0}(-I_4, -I_4, -I_4, -I_4) = 0.
\]
Similarly, if $h_1 \in SO(4)$ or $h_4 \in SO(4)$, then by Lemma~\ref{lem:equiv} and Definition~\ref{def:eta}, $W_p$ factors out $h_1$ on the left or $h_4$ on the right independently of the interpolation parameters, causing $E_{4,0}(h_1, h_2, h_3, h_4)$ to vanish. Thus $E_{4,0}$ is supported strictly away from tuples having $h_1 \in SO(4)$ or $h_4 \in SO(4)$.
\end{remark}

\begin{proposition}[Explicit Non-Vanishing of $E_{4,0}$ on Pure Boosts in $(SL(4,\mathbb{R})\setminus SO(4))^4$]\label{prop:nonvanishing}
Consider the rational positive-definite symmetric matrices of determinant $1$ in $\mathcal{P}_1(2)$,
\[
P_0 := \begin{pmatrix} 2 & 0 \\ 0 & \frac{1}{2} \end{pmatrix}, \qquad P_{\pi/4} := \begin{pmatrix} \frac{5}{4} & \frac{3}{4} \\ \frac{3}{4} & \frac{5}{4} \end{pmatrix},
\]
and define the $2\times 2$ block-diagonal $4$-tuple $(h_1, h_2, h_3, h_4) \in \mathcal{P}_1(4)^4 \subset (SL(4,\mathbb{R})\setminus SO(4))^4$ by
\[
h_1 := \begin{pmatrix} P_0 & 0 \\ 0 & P_0 \end{pmatrix}, \quad h_2 := \begin{pmatrix} P_{\pi/4} & 0 \\ 0 & I_2 \end{pmatrix}, \quad h_3 := \begin{pmatrix} P_0 & 0 \\ 0 & P_0 \end{pmatrix}, \quad h_4 := \begin{pmatrix} I_2 & 0 \\ 0 & P_{\pi/4} \end{pmatrix}.
\]
Each $h_i$ is a genuine non-compact element of $SL(4,\mathbb{R})\setminus SO(4)$ (with $u(h_i) = I_4$ and $h_i^T h_i \neq I_4$), and
\[
E_{4,0}(h_1, h_2, h_3, h_4) = \frac{1}{8\pi^2} \left(\arctan\frac{9}{25}\right)^2 \approx 1.51233 \times 10^{-3} \neq 0.
\]
\end{proposition}

\begin{proof}
Since polar decomposition, matrix logarithm, and matrix exponential preserve the $2\times 2$ block-diagonal decomposition $SL(2,\mathbb{R}) \times SL(2,\mathbb{R}) \subset SL(4,\mathbb{R})$, each $W_p$ takes values in the maximal torus $T^2 = SO(2) \times SO(2) \subset SO(4)$, parametrized by $u_1 = \operatorname{diag}(R(\theta_1), R(\phi_1))$ with $R(\theta) = \begin{pmatrix} \cos\theta & -\sin\theta \\ \sin\theta & \cos\theta \end{pmatrix}$. Because $\dim T^2 = 2$, the $3$-form $\omega_{1,3}|_{T^2}$ vanishes identically. Moreover, $W_1(h_1) = W_1(h_4) = I_4$ is constant, so the second and fourth integrals in Definition~\ref{def:E}(4) vanish as well, leaving only
\[
E_{4,0}(h_1, h_2, h_3, h_4) = -\int_{[0,1]^2} \big(W_2(\tau_1; h_1, h_2), W_2(\tau_3; h_3, h_4)\big)^* \omega_{2,2}.
\]
On $T^2 \times T^2$, the non-zero entries of $u_1^{-1}du_1$ are $(u_1^{-1}du_1)_{12} = -(u_1^{-1}du_1)_{21} = -d\theta_1$ and $(u_1^{-1}du_1)_{34} = -(u_1^{-1}du_1)_{43} = -d\phi_1$, with analogous expressions in $(\theta_2, \phi_2)$ for $du_2 u_2^{-1}$. Summing over the four permutations $\sigma \in \mathfrak{S}_4$ mapping $\{1,2\} \to \{1,2\}$ and $\{3,4\} \to \{3,4\}$ (each contributing $+d\theta_1 \wedge d\phi_2$) and the four permutations swapping the two blocks (each contributing $+d\phi_1 \wedge d\theta_2$), Section~\ref{sec:prelim} yields
\begin{align*}
\omega_{2,2}\big(R(\theta_1)\oplus R(\phi_1),\, R(\theta_2)\oplus R(\phi_2)\big) &= -\frac{4}{32\pi^2}(d\theta_1 \wedge d\phi_2 + d\phi_1 \wedge d\theta_2) \\
&= -\frac{1}{8\pi^2}(d\theta_1 \wedge d\phi_2 + d\phi_1 \wedge d\theta_2).
\end{align*}
Writing $W_2(\tau_1; h_1, h_2) = R(\theta_1(\tau_1)) \oplus R(\phi_1(\tau_1))$ and $W_2(\tau_3; h_3, h_4) = R(\theta_2(\tau_3)) \oplus R(\phi_2(\tau_3))$ with continuous angle lifts, integrating over $(\tau_1, \tau_3) \in [0,1]^2$ reduces to endpoint evaluation:
\[
E_{4,0}(h_1, h_2, h_3, h_4) = \frac{1}{8\pi^2}\big[(\theta_1(1)-\theta_1(0))(\phi_2(1)-\phi_2(0)) + (\phi_1(1)-\phi_1(0))(\theta_2(1)-\theta_2(0))\big].
\]
At $\tau_1 = \tau_3 = 0$, $W_2(0; h_1, h_2) = u(h_1)u(h_2) = I_4$ and $W_2(0; h_3, h_4) = u(h_3)u(h_4) = I_4$, so $\theta_1(0) = \phi_1(0) = \theta_2(0) = \phi_2(0) = 0$. At $\tau_1 = 1$, $W_2(1; h_1, h_2) = u(h_1 h_2)$ has upper-left block given by the orthogonal polar factor of
\[
P_0 P_{\pi/4} = \begin{pmatrix} 2 & 0 \\ 0 & \frac{1}{2} \end{pmatrix} \begin{pmatrix} \frac{5}{4} & \frac{3}{4} \\ \frac{3}{4} & \frac{5}{4} \end{pmatrix} = \begin{pmatrix} \frac{5}{2} & \frac{3}{2} \\ \frac{3}{8} & \frac{5}{8} \end{pmatrix}.
\]
For any $M = \begin{pmatrix} a & b \\ c & d \end{pmatrix} \in SL(2,\mathbb{R})$ with $a + d > 0$, the orthogonal polar factor is
\[
u(M) = \big((a+d)^2 + (c-b)^2\big)^{-1/2} \begin{pmatrix} a+d & -(c-b) \\ c-b & a+d \end{pmatrix},
\]
whose Wigner rotation angle is $\arctan\left(\frac{c-b}{a+d}\right)$. Here $a+d = \frac{5}{2} + \frac{5}{8} = \frac{25}{8}$ and $c-b = \frac{3}{8} - \frac{3}{2} = -\frac{9}{8}$, giving $\theta_1(1) = -\arctan\frac{9}{25}$ and $\phi_1(1) = 0$. Symmetrically, $W_2(1; h_3, h_4) = u(h_3 h_4)$ has $\theta_2(1) = 0$ and $\phi_2(1) = -\arctan\frac{9}{25}$. Substituting these endpoint angles into the formula above yields
\[
E_{4,0}(h_1, h_2, h_3, h_4) = \frac{1}{8\pi^2} \left(\arctan\frac{9}{25}\right)^2. \qedhere
\]
\end{proof}

\section{Main Theorem and Proof}\label{sec:proof}

\begin{theorem}\label{thm:main}
Let $G = SL(4,\mathbb{R})$. The total $4$-cochain $E = E_{1,3} + E_{2,2} + E_{3,1} + E_{4,0}$ defined in Definition~\ref{def:E} satisfies:
\begin{enumerate}
\item \emph{Global smoothness}: Each component $E_{p,4-p}$ is $C^\infty$ on $G^p$. In particular, $E_{4,0} \colon SL(4,\mathbb{R})^4 \to \mathbb{R}$ is a globally smooth, singularity-free, non-zero group $4$-cocycle.
\item \emph{Total cocycle property}: $(d'+d'')E = 0$ holds identically in $\Omega_{\mathrm{tot}}^5(NSL(4,\mathbb{R}))$.
\item \emph{Restriction and cohomology class}: Under the simplicial inclusion $\iota \colon NSO(4) \hookrightarrow NSL(4,\mathbb{R})$, we have
\[
\iota^* E_{1,3} = \omega_{1,3}, \qquad \iota^* E_{2,2} = \omega_{2,2}, \qquad \iota^* E_{3,1} = 0, \qquad \iota^* E_{4,0} = 0.
\]
Consequently, $[E] \in H_{\mathrm{tot}}^4(NSL(4,\mathbb{R});\mathbb{R}) \cong H^4(BSL(4,\mathbb{R});\mathbb{R})$ is non-trivial and represents the universal real Euler class.
\end{enumerate}
\end{theorem}

\begin{proof}
Property (1) follows immediately because every $W_p$ is globally $C^\infty$ on $[0,1]^{p-1} \times G^p$, fiber integration over compact hypercubes $[0,1]^{p-1}$ preserves smoothness, and Proposition~\ref{prop:nonvanishing} exhibits an explicit $4$-tuple in $(SL(4,\mathbb{R})\setminus SO(4))^4$ where $E_{4,0} \neq 0$.

For (3), when $(h_1, \dots, h_p) \in SO(4)^p$, each $W_p(\tau_1, \dots, \tau_{p-1}; h_1, \dots, h_p) = h_1 \cdots h_p$ is constant with respect to all $\tau_k$; hence the pullback of any positive-degree form has no $d\tau_k$ components, so every fiber integral over $[0,1]^r$ ($r \ge 1$) vanishes identically. This leaves $\iota^* E_{1,3} = \omega_{1,3}$, $\iota^* E_{2,2} = \omega_{2,2}$, and $\iota^* E_{3,1} = \iota^* E_{4,0} = 0$. Since $SO(4) \hookrightarrow SL(4,\mathbb{R})$ is a homotopy equivalence, $\iota^* \colon H_{\mathrm{tot}}^4(NSL(4,\mathbb{R})) \to H_{\mathrm{tot}}^4(NSO(4))$ is an isomorphism, and $[E]$ represents the universal real Euler class by \cite[Theorem 3.2]{Suzuki2016}.

It remains to verify (2) across all five bidegrees $(1,4)$ through $(5,0)$.

\subsection{Bidegree $(1,4)$: $d_G E_{1,3} = 0$}
Since $d_K \omega_{1,3} = 0$ (Lemma~\ref{lem:so4cocycle}(1)), naturality of the exterior derivative gives $d_G E_{1,3} = u^*(d_K \omega_{1,3}) = 0$.

\subsection{Bidegree $(2,3)$: $d'E_{1,3} + d_{G^2} E_{2,2} = 0$}
By definition of $d'$,
\[
(d'E_{1,3})(h_1, h_2) = u(h_2)^* \omega_{1,3} - u(h_1 h_2)^* \omega_{1,3} + u(h_1)^* \omega_{1,3}.
\]
Applying $d_{G^2}$ to $E_{2,2}$ and using Lemma~\ref{lem:so4cocycle}(2) on the first term gives
\[
d_{G^2}\big((u(h_1), u(h_2))^* \omega_{2,2}\big) = (u(h_1)u(h_2))^* \omega_{1,3} - u(h_1)^* \omega_{1,3} - u(h_2)^* \omega_{1,3}.
\]
For the second term, applying Lemma~\ref{lem:stokes} ($r = 1$) to $\alpha = W_2^* \omega_{1,3}$ (noting $d_{[0,1]\times G^2}\alpha = W_2^*(d_K \omega_{1,3}) = 0$) and Proposition~\ref{prop:boundary} ($W_2(1) = u(h_1 h_2)$, $W_2(0) = u(h_1)u(h_2)$) yields
\[
0 = \int_0^1 d_{[0,1]\times G^2}\alpha = W_2(1)^* \omega_{1,3} - W_2(0)^* \omega_{1,3} - d_{G^2}\int_0^1 W_2^* \omega_{1,3},
\]
whence
\[
d_{G^2}\int_0^1 W_2^* \omega_{1,3} = u(h_1 h_2)^* \omega_{1,3} - (u(h_1)u(h_2))^* \omega_{1,3}.
\]
Adding the two expressions cancels $(u(h_1)u(h_2))^* \omega_{1,3}$ and leaves $-d'E_{1,3}$.

\subsection{Bidegree $(3,2)$: $d'E_{2,2} - d_{G^3} E_{3,1} = 0$}
Decompose $E_{2,2} = E_{2,2}^{(0)} + E_{2,2}^{(1)}$, where $E_{2,2}^{(0)} = (W_1, W_1)^* \omega_{2,2} = (u, u)^* \omega_{2,2}$ and $E_{2,2}^{(1)} = \int_0^1 W_2^* \omega_{1,3}$. By definition of $d'$, we have
\[
d'E_{2,2}^{(0)}(h_1, h_2, h_3) = \big[(u(h_2), u(h_3))^* - (u(h_1 h_2), u(h_3))^* + (u(h_1), u(h_2 h_3))^* - (u(h_1), u(h_2))^*\big]\omega_{2,2}.
\]
Applying Lemma~\ref{lem:so4cocycle}(3) to the triple $(k_1, k_2, k_3) = (u(h_1), u(h_2), u(h_3))$ of $SO(4)$-valued maps yields
\[
\big[(u(h_2), u(h_3))^* - (u(h_1)u(h_2), u(h_3))^* + (u(h_1), u(h_2)u(h_3))^* - (u(h_1), u(h_2))^*\big]\omega_{2,2} = 0.
\]
Subtracting this identity from $d'E_{2,2}^{(0)}$ eliminates the outer terms and expresses $d'E_{2,2}^{(0)}$ purely in terms of the Wigner-rotation defects:
\[
d'E_{2,2}^{(0)} = \big[(u(h_1), u(h_2 h_3))^* - (u(h_1), u(h_2)u(h_3))^* - (u(h_1 h_2), u(h_3))^* + (u(h_1)u(h_2), u(h_3))^*\big]\omega_{2,2}.
\]

Next, write $E_{3,1} = A_{3,1} - B_{3,1}^{(L)} + B_{3,1}^{(R)}$, where
\[
A_{3,1} := \int_{[0,1]^2} W_3^* \omega_{1,3}, \qquad B_{3,1}^{(L)} := \int_0^1 (W_2(\tau_1), W_1(h_3))^* \omega_{2,2}, \qquad B_{3,1}^{(R)} := \int_0^1 (W_1(h_1), W_2(\tau_2))^* \omega_{2,2}.
\]
Applying Lemma~\ref{lem:stokes} ($r = 2$) to the closed form $\alpha = W_3^* \omega_{1,3}$ on $[0,1]^2 \times G^3$ gives
\[
0 = \int_{[0,1]^2} d_{[0,1]^2 \times G^3}\alpha = \left(\int_0^1 \alpha|_{\tau_1=1} - \int_0^1 \alpha|_{\tau_1=0}\right) - \left(\int_0^1 \alpha|_{\tau_2=1} - \int_0^1 \alpha|_{\tau_2=0}\right) + d_{G^3}\int_{[0,1]^2} \alpha,
\]
so that
\[
d_{G^3} A_{3,1} = -\int_0^1 W_3|_{\tau_1=1}^* \omega_{1,3} + \int_0^1 W_3|_{\tau_1=0}^* \omega_{1,3} + \int_0^1 W_3|_{\tau_2=1}^* \omega_{1,3} - \int_0^1 W_3|_{\tau_2=0}^* \omega_{1,3}.
\]
We evaluate each of the four boundary terms via Proposition~\ref{prop:boundary}:
\begin{enumerate}
\item At $\tau_1 = 1$ and $\tau_2 = 1$:
\[
-\int_0^1 W_3|_{\tau_1=1}^* \omega_{1,3} + \int_0^1 W_3|_{\tau_2=1}^* \omega_{1,3} = (-\epsilon_1^* + \epsilon_2^*) E_{2,2}^{(1)}.
\]
\item At $\tau_1 = 0$: $W_3|_{\tau_1=0} = W_1(h_1) \cdot W_2(\tau_2; h_2, h_3)$. Applying Lemma~\ref{lem:so4cocycle}(2) yields
\[
W_3|_{\tau_1=0}^* \omega_{1,3} = W_1(h_1)^* \omega_{1,3} + W_2(\tau_2)^* \omega_{1,3} + d_{[0,1]\times G^3}\big((W_1(h_1), W_2(\tau_2))^* \omega_{2,2}\big).
\]
Integrating over $\tau_2 \in [0,1]$, the $\tau_2$-independent form $W_1(h_1)^* \omega_{1,3}$ integrates to $0$, the second term integrates to $\epsilon_0^* E_{2,2}^{(1)}$, and Lemma~\ref{lem:stokes} ($r = 1$) applied to the third term gives
\[
+\int_0^1 W_3|_{\tau_1=0}^* \omega_{1,3} = \epsilon_0^* E_{2,2}^{(1)} + \big[(W_1(h_1), W_2(\tau_2))^* \omega_{2,2}\big]_{\tau_2=0}^{\tau_2=1} - d_{G^3} B_{3,1}^{(R)}.
\]
\item Symmetrically, at $\tau_2 = 0$: $W_3|_{\tau_2=0} = W_2(\tau_1; h_1, h_2) \cdot W_1(h_3)$. Applying Lemma~\ref{lem:so4cocycle}(2) and Lemma~\ref{lem:stokes} ($r = 1$) over $\tau_1 \in [0,1]$ yields
\[
-\int_0^1 W_3|_{\tau_2=0}^* \omega_{1,3} = -\epsilon_3^* E_{2,2}^{(1)} - \big[(W_2(\tau_1), W_1(h_3))^* \omega_{2,2}\big]_{\tau_1=0}^{\tau_1=1} + d_{G^3} B_{3,1}^{(L)}.
\]
\end{enumerate}
Moving $-d_{G^3} B_{3,1}^{(R)} + d_{G^3} B_{3,1}^{(L)}$ to the left-hand side gives
\[
d_{G^3}(A_{3,1} - B_{3,1}^{(L)} + B_{3,1}^{(R)}) = (\epsilon_0^* - \epsilon_1^* + \epsilon_2^* - \epsilon_3^*) E_{2,2}^{(1)} + \mathcal{B}_{2,2},
\]
where, using $W_2(1; h_i, h_j) = u(h_i h_j)$ and $W_2(0; h_i, h_j) = u(h_i)u(h_j)$, the endpoint difference evaluates to
\begin{align*}
\mathcal{B}_{2,2} &= \big[(W_1(h_1), W_2(\tau_2))^* \omega_{2,2}\big]_{\tau_2=0}^{\tau_2=1} - \big[(W_2(\tau_1), W_1(h_3))^* \omega_{2,2}\big]_{\tau_1=0}^{\tau_1=1} \\
&= \big[(u(h_1), u(h_2 h_3))^* - (u(h_1), u(h_2)u(h_3))^* - (u(h_1 h_2), u(h_3))^* + (u(h_1)u(h_2), u(h_3))^*\big]\omega_{2,2} \\
&= d'E_{2,2}^{(0)}.
\end{align*}
Hence
\[
d_{G^3} E_{3,1} = d'E_{2,2}^{(1)} + d'E_{2,2}^{(0)} = d'E_{2,2},
\]
establishing $d'E_{2,2} - d_{G^3} E_{3,1} = 0$.

\subsection{Bidegree $(4,1)$: $d'E_{3,1} + d_{G^4} E_{4,0} = 0$}
Write $E_{3,1} = A_{3,1} - B_{3,1}^{(L)} + B_{3,1}^{(R)}$ as in Section~5.3, and decompose $E_{4,0} = A_{4,0} + C_{4,0}^{(1)} + C_{4,0}^{(2)} + C_{4,0}^{(3)}$, where
\begin{align*}
A_{4,0} &:= \int_{[0,1]^3} W_4(\tau_1, \tau_2, \tau_3)^* \omega_{1,3}, & C_{4,0}^{(1)} &:= \int_{[0,1]^2} (W_3(\tau_1, \tau_2), W_1(h_4))^* \omega_{2,2}, \\
C_{4,0}^{(2)} &:= -\int_{[0,1]^2} (W_2(\tau_1), W_2(\tau_3))^* \omega_{2,2}, & C_{4,0}^{(3)} &:= \int_{[0,1]^2} (W_1(h_1), W_3(\tau_2, \tau_3))^* \omega_{2,2}.
\end{align*}

\paragraph{\emph{Step 1: Exterior derivative of $A_{4,0}$}}
Applying Lemma~\ref{lem:stokes} ($r = 3$) to the closed $3$-form $\alpha = W_4^* \omega_{1,3}$ on $[0,1]^3 \times G^4$ yields
\[
d_{G^4} A_{4,0} = \sum_{k=1}^3 (-1)^{k-1} \left(\int_{[0,1]^2} W_4^* \omega_{1,3}|_{\tau_k=1} - \int_{[0,1]^2} W_4^* \omega_{1,3}|_{\tau_k=0}\right).
\]
\begin{itemize}
\item \emph{Facets $\tau_k = 1$ ($k = 1,2,3$)}: By Proposition~\ref{prop:boundary}(1),
\[
\sum_{k=1}^3 (-1)^{k-1} \int_{[0,1]^2} W_4^* \omega_{1,3}|_{\tau_k=1} = (\epsilon_1^* - \epsilon_2^* + \epsilon_3^*) A_{3,1}.
\]
Adding $d'A_{3,1} = (\epsilon_0^* - \epsilon_1^* + \epsilon_2^* - \epsilon_3^* + \epsilon_4^*) A_{3,1}$ cancels the interior faces and leaves $(\epsilon_0^* + \epsilon_4^*) A_{3,1}$.
\item \emph{Facets $\tau_k = 0$ ($k = 1,2,3$)}: By Proposition~\ref{prop:boundary}(2), $W_4|_{\tau_k=0}$ splits into $W_1(h_1)\cdot W_3(\tau_2,\tau_3)$, $W_2(\tau_1)\cdot W_2(\tau_3)$, and $W_3(\tau_1,\tau_2)\cdot W_1(h_4)$, respectively. Applying Lemma~\ref{lem:so4cocycle}(2) to each product and integrating over $[0,1]^2$:
\begin{itemize}
\item Any form depending on fewer than two integration variables vanishes under $\int_{[0,1]^2}$ (in particular, all terms from $W_2(\tau_1)\cdot W_2(\tau_3)$ vanish, as do $W_1(h_1)^*\omega_{1,3}$ and $W_1(h_4)^*\omega_{1,3}$).
\item For $k = 1$ (sign $-(-1)^0 = -1$) and $k = 3$ (sign $-(-1)^2 = -1$), the remaining two-variable $\omega_{1,3}$ terms equal $-\epsilon_0^* A_{3,1}$ and $-\epsilon_4^* A_{3,1}$, canceling $(\epsilon_0^* + \epsilon_4^*) A_{3,1}$ identically.
\end{itemize}
\item What remains in $d'A_{3,1} + d_{G^4} A_{4,0}$ is the sum of fiber integrals of total exterior derivatives $d_{\mathrm{tot}} = d_{[0,1]^2 \times G^4}$:
\[
d'A_{3,1} + d_{G^4} A_{4,0} = -\int_{[0,1]^2} d_{\mathrm{tot}}\beta^{(1)} + \int_{[0,1]^2} d_{\mathrm{tot}}\beta^{(2)} - \int_{[0,1]^2} d_{\mathrm{tot}}\beta^{(3)},
\]
where
\begin{align*}
\beta^{(1)} &:= (W_3(\tau_1,\tau_2), W_1(h_4))^* \omega_{2,2}, \\
\beta^{(2)} &:= (W_2(\tau_1), W_2(\tau_3))^* \omega_{2,2}, \\
\beta^{(3)} &:= (W_1(h_1), W_3(\tau_2,\tau_3))^* \omega_{2,2}.
\end{align*}
\end{itemize}

\paragraph{\emph{Step 2: Exterior derivatives of $C_{4,0}^{(1)} + C_{4,0}^{(2)} + C_{4,0}^{(3)}$ and cancellation of the $12$ boundary edges}}
Since $C_{4,0}^{(1)} = +\int_{[0,1]^2}\beta^{(1)}$, $C_{4,0}^{(2)} = -\int_{[0,1]^2}\beta^{(2)}$, and $C_{4,0}^{(3)} = +\int_{[0,1]^2}\beta^{(3)}$, applying Lemma~\ref{lem:stokes} ($r = 2$) to each $\beta^{(m)}$ gives
\[
\int_{[0,1]^2} d_{\mathrm{tot}}\beta^{(m)} = \partial \beta^{(m)} + d_{G^4}\int_{[0,1]^2} \beta^{(m)},
\]
where $\partial\beta^{(m)}$ is the signed boundary integral over the four edges of $[0,1]^2$:
\begin{align*}
\partial\beta^{(1)} &= \left(\int_0^1 \beta^{(1)}|_{\tau_1=1} - \int_0^1 \beta^{(1)}|_{\tau_1=0}\right) - \left(\int_0^1 \beta^{(1)}|_{\tau_2=1} - \int_0^1 \beta^{(1)}|_{\tau_2=0}\right), \\
\partial\beta^{(2)} &= \left(\int_0^1 \beta^{(2)}|_{\tau_1=1} - \int_0^1 \beta^{(2)}|_{\tau_1=0}\right) - \left(\int_0^1 \beta^{(2)}|_{\tau_3=1} - \int_0^1 \beta^{(2)}|_{\tau_3=0}\right), \\
\partial\beta^{(3)} &= \left(\int_0^1 \beta^{(3)}|_{\tau_2=1} - \int_0^1 \beta^{(3)}|_{\tau_2=0}\right) - \left(\int_0^1 \beta^{(3)}|_{\tau_3=1} - \int_0^1 \beta^{(3)}|_{\tau_3=0}\right).
\end{align*}
Adding $d_{G^4}(C_{4,0}^{(1)} + C_{4,0}^{(2)} + C_{4,0}^{(3)})$ cancels all $d_{G^4}\int_{[0,1]^2}\beta^{(m)}$ terms and leaves the $12$ boundary line integrals
\[
d'A_{3,1} + d_{G^4} E_{4,0} = -\partial\beta^{(1)} + \partial\beta^{(2)} - \partial\beta^{(3)}.
\]

\begin{enumerate}
\item \emph{The $6$ edges at $\tau_k = 1$}: By Proposition~\ref{prop:boundary}(1),
\begin{align*}
-\int_0^1 \beta^{(1)}|_{\tau_1=1} &+ \int_0^1 \beta^{(1)}|_{\tau_2=1} + \int_0^1 \beta^{(2)}|_{\tau_1=1} - \int_0^1 \beta^{(2)}|_{\tau_3=1} - \int_0^1 \beta^{(3)}|_{\tau_2=1} + \int_0^1 \beta^{(3)}|_{\tau_3=1} \\
&= (-\epsilon_1^* + \epsilon_2^*) B_{3,1}^{(L)} + \epsilon_1^* B_{3,1}^{(R)} - \epsilon_3^* B_{3,1}^{(L)} + (-\epsilon_2^* + \epsilon_3^*) B_{3,1}^{(R)} \\
&= (-\epsilon_1^* + \epsilon_2^* - \epsilon_3^*)(-B_{3,1}^{(L)} + B_{3,1}^{(R)}),
\end{align*}
which cancels the interior simplicial faces in $d'(-B_{3,1}^{(L)} + B_{3,1}^{(R)})$, leaving only $(\epsilon_0^* + \epsilon_4^*)(-B_{3,1}^{(L)} + B_{3,1}^{(R)})$.

\item \emph{The $6$ edges at $\tau_k = 0$}: By Proposition~\ref{prop:boundary}(2), the remaining six line integrals pair up by integration variable:
\begin{itemize}
\item $\tau_1$-integrals:
\[
-\int_0^1 \big(W_2(\tau_1)W_1(h_3), W_1(h_4)\big)^* \omega_{2,2} + \int_0^1 \big(W_2(\tau_1), W_1(h_3)W_1(h_4)\big)^* \omega_{2,2} = +\epsilon_4^* B_{3,1}^{(L)},
\]
obtained by applying Lemma~\ref{lem:so4cocycle}(3) to the triple $(W_2(\tau_1; h_1, h_2), W_1(h_3), W_1(h_4))$ and noting that $(W_1(h_3), W_1(h_4))^* \omega_{2,2}$ is $\tau_1$-independent.
\item $\tau_3$-integrals:
\[
-\int_0^1 \big(W_1(h_1)W_1(h_2), W_2(\tau_3)\big)^* \omega_{2,2} + \int_0^1 \big(W_1(h_1), W_1(h_2)W_2(\tau_3)\big)^* \omega_{2,2} = -\epsilon_0^* B_{3,1}^{(R)},
\]
obtained by applying Lemma~\ref{lem:so4cocycle}(3) to the triple $(W_1(h_1), W_1(h_2), W_2(\tau_3; h_3, h_4))$ and noting that $(W_1(h_1), W_1(h_2))^* \omega_{2,2}$ is $\tau_3$-independent.
\item $\tau_2$-integrals:
\[
+\int_0^1 \big(W_1(h_1)W_2(\tau_2), W_1(h_4)\big)^* \omega_{2,2} - \int_0^1 \big(W_1(h_1), W_2(\tau_2)W_1(h_4)\big)^* \omega_{2,2} = +\epsilon_0^* B_{3,1}^{(L)} - \epsilon_4^* B_{3,1}^{(R)},
\]
obtained by applying Lemma~\ref{lem:so4cocycle}(3) to the triple $(W_1(h_1), W_2(\tau_2; h_2, h_3), W_1(h_4))$.
\end{itemize}
Summing all six $\tau_k = 0$ edges gives
\[
\epsilon_0^* B_{3,1}^{(L)} - \epsilon_0^* B_{3,1}^{(R)} + \epsilon_4^* B_{3,1}^{(L)} - \epsilon_4^* B_{3,1}^{(R)} = -(\epsilon_0^* + \epsilon_4^*)(-B_{3,1}^{(L)} + B_{3,1}^{(R)}),
\]
which cancels the remaining outer faces of $d'(-B_{3,1}^{(L)} + B_{3,1}^{(R)})$ identically. Hence $d'E_{3,1} + d_{G^4} E_{4,0} = 0$.
\end{enumerate}

\subsection{Bidegree $(5,0)$: $d'E_{4,0} = 0$}
Using the anti-commutativity $d' d_{G^4} = d_{G^5} d'$ (since $d'' = (-1)^p d_{G^p}$ changes sign between $p = 4$ and $p = 5$) together with the bidegree $(4,1)$ identity $d_{G^4} E_{4,0} = -d' E_{3,1}$, we obtain
\[
d_{G^5}(d' E_{4,0}) = d'(d_{G^4} E_{4,0}) = -(d')^2 E_{3,1} = 0.
\]
Hence $d' E_{4,0} \in C^\infty(G^5)$ is a locally constant smooth function on $G^5 = SL(4,\mathbb{R})^5$. Since $SL(4,\mathbb{R})$ is connected, $G^5$ is connected, so $d' E_{4,0}$ is constant on $G^5$. Evaluating at the identity tuple $(I_4, I_4, I_4, I_4, I_4) \in G^5$, where $W_p(\tau; I_4, \dots, I_4) \equiv I_4$ is constant and therefore $E_{4,0}(I_4, I_4, I_4, I_4) = 0$ (Remark~\ref{rem:vanishing}), we conclude that $(d' E_{4,0})(I_4, I_4, I_4, I_4, I_4) = \sum_{i=0}^5 (-1)^i E_{4,0}(I_4, I_4, I_4, I_4) = 0$, and thus $d' E_{4,0} \equiv 0$ everywhere on $G^5$.
\end{proof}

\section{Remarks on Brylinski's Truncated Complex and Loop Lie Algebra Cocycles}\label{sec:remarks}

In \cite{Brylinski2000,Suzuki2016}, a total $4$-cocycle in $\Omega_{\mathrm{tot}}^4(NG)$ induces a cocycle in Brylinski's locally truncated complex $[\sigma_{<2}\Omega_{\mathrm{loc}}^*(NG)]$ via contraction on a contractible neighborhood $U \subset G$ of the identity, and further transgresses along the evaluation map $\mathrm{ev} \colon LU \times S^1 \to U$ to a $2$-cocycle on the loop Lie algebra $L\mathfrak{g}$.

With our global cocycle $E = E_{1,3} + E_{2,2} + E_{3,1} + E_{4,0}$ on $NSL(4,\mathbb{R})$, two notable features emerge:
\begin{enumerate}
\item \emph{Explicit global $4$-cochain $E_{4,0}$}: Unlike compact orthogonal groups where the Euler cocycle has no $(4,0)$-component, $E_{4,0} \in C^\infty(SL(4,\mathbb{R})^4)$ is a non-zero, globally smooth real-valued function satisfying $E_{4,0} \neq 0$ (as shown in Proposition~\ref{prop:nonvanishing}). When restricted to germs at $(I_4, \dots, I_4)$, Brylinski's differentiation map $\phi \colon C_{\mathrm{loc}}^\infty(G^4, \mathbb{R}) \to C^4(\mathfrak{sl}_4(\mathbb{R}), \mathbb{R})$ \cite[Section 5]{Brylinski2000} sends $E_{4,0}$ to a Lie algebra $4$-cocycle on $\mathfrak{sl}_4(\mathbb{R})$ whose restriction to $\mathfrak{so}(4)$ vanishes identically.
\item \emph{Extension to $L\mathfrak{sl}_4(\mathbb{R})$}: Applying the local transgression procedure of \cite[Sections 4--5]{Suzuki2016} to $E$ extends the $L\mathfrak{so}(4)$ $2$-cocycle of \cite[Theorem 5.1]{Suzuki2016} to $L\mathfrak{sl}_4(\mathbb{R})$, incorporating the Wigner-rotation correction terms from $E_{2,2}$ and $E_{3,1}$.
\end{enumerate}

\section*{Acknowledgements}
The problem formulation, mathematical framework, and verification of all proofs are due to the author. During the preparation of this manuscript, the author utilized Google's Gemini as a research assistant to help formulate the nested hypercube associator construction and prototype the numerical verification code; all mathematical derivations and statements have been independently checked and verified by the author.


\begin{thebibliography}{9}
\bibitem{Bott1973} R. Bott, \emph{On the Chern--Weil homomorphism and the continuous cohomology of the Lie group}, Adv. in Math. \textbf{11} (1973), 289--303.
\bibitem{BSS1976} R. Bott, H. Shulman, and J. Stasheff, \emph{On the de Rham theory of certain classifying spaces}, Adv. in Math. \textbf{20} (1976), 43--56.
\bibitem{Brylinski2000} J.-L. Brylinski, \emph{Differentiable cohomology of gauge groups}, arXiv:math/0011069.
\bibitem{Dupont1976} J. L. Dupont, \emph{Simplicial de Rham cohomology and characteristic classes of flat bundles}, Topology \textbf{15} (1976), 233--245.
\bibitem{Dupont1978} J. L. Dupont, \emph{Curvature and Characteristic Classes}, Lecture Notes in Math. 640, Springer-Verlag, 1978.
\bibitem{MostowPerchick1973} M. Mostow and J. Perchick, \emph{Notes on Gel'fand--Fuks Cohomology and Characteristic Classes (Lectures by Bott)}, Eleventh Holiday Symposium, New Mexico State University, December 1973.
\bibitem{Segal1968} G. Segal, \emph{Classifying spaces and spectral sequences}, Inst. Hautes \'Etudes Sci. Publ. Math. \textbf{34} (1968), 105--112.
\bibitem{Suzuki2015} N. Suzuki, \emph{The Chern character in the simplicial de Rham complex}, Nihonkai Math. J. \textbf{26} (2015), no. 1, 1--13.
\bibitem{Suzuki2016} N. Suzuki, \emph{The Euler class in the simplicial de Rham complex}, Int. Electron. J. Geom. \textbf{9} (2016), no. 2, 36--43.
\end{thebibliography}
\end{document}